\documentclass[11pt]{amsart}
\usepackage[margin=1.4in]{geometry}
\usepackage{amsmath,amsfonts,amssymb,amsthm,mathtools}
\usepackage{cases}
\usepackage{url}
\usepackage[all]{xy}
\usepackage{tikz-cd}
\usepackage{array}
\usepackage{mathtools}
\usepackage{faktor}
\usepackage{caption} 
\usepackage{verbatim,eucal}
\usepackage{xfrac}
\usepackage{MnSymbol}

\usepackage{float}
\usepackage{multicol}

\usepackage{letltxmacro}
\LetLtxMacro\Oldfootnote\footnote

\usepackage{soul}

\usepackage{hyperref}
\hypersetup{colorlinks=true,citecolor=blue!70!black,linkcolor=red!60!black,linktocpage=true}
\usepackage[capitalise,nameinlink]{cleveref}
\crefname{enumi}{part}{parts}

\makeatother
\numberwithin{equation}{section}

\newtheorem{thm}[equation]{Theorem} 
\newtheorem{prop}[equation]{Proposition}
\crefname{prop}{Proposition}{Propositions}

\newtheorem{lemma}[equation]{Lemma} 
\newtheorem{cor}[equation]{Corollary}
\crefname{cor}{Corollary}{Corollaries}

\newtheorem{example}[equation]{Example}
\newtheorem{remark}[equation]{Remark}

\newenvironment{rem}{\begin{remark}\rm}{\end{remark}}

\DeclareMathOperator{\Ext}{Ext}
\DeclareMathOperator{\gr}{gr} 

\DeclareMathOperator{\im}{Im}
\DeclareMathOperator{\Ker}{Ker}

\DeclareMathOperator{\Span}{Span}

\newcommand{\NN}{\mathbb N}
\newcommand{\FF}{\mathbb F}

\newcommand{\DOT}{\setlength{\unitlength}{1pt}\begin{picture}(2.5,2)
          (1,1)\put(2,3.5){\circle*{3}}\end{picture}}

\newcommand{\id}{\mbox{\rm id}}
\newcommand{\Hom}{\mbox{\rm Hom\,}}

\renewcommand{\ker}{\mbox{\rm Ker\,}}

\newcommand{\ot}{\otimes}

\newcommand{\cH}{\mathcal{H}}

\DeclareMathOperator{\codim}{codim}

\newcommand{\HH}{{\rm HH}}
\newcommand{\Wedge}{
{\textstyle \bigwedge}}

\newcommand{\ld}{\lambda}   

\definecolor{SFA}{HTML}{CB00F5}

\begin{document}


\begin{abstract}
We investigate deformations of skew group algebras that arise from a finite cyclic group acting on a polynomial ring in positive characteristic, where characteristic divides the order of the group.  We allow deformations which deform both the group action and the vector space multiplication.  We fully characterize the Poincar\'{e}-Birkhoff-Witt deformations which arise in this setting from multiple perspectives: a necessary and sufficient condition list, a practical road map from which one can generate examples corresponding to any choice of group algebra element, an explicit formula, and a combinatorial analysis of the class of algebras.
\end{abstract}

\title[Graded deformations of skew group algebras]{Graded deformations of skew group algebras for cyclic\\
transvection
groups acting on polynomial rings 
\\
in positive characteristic}

\date{\today}

\author{Lauren Grimley}
\address{Lauren Grimley, Department of Mathematics and School of Computer Science, University of Oklahoma, Norman, OK}
\email{lgrimley@ou.edu}

\author{Naomi Krawzik}
\address{Naomi Krawzik, Department of Mathematics and Statistics, Sam Houston State University,
Huntsville, Texas, USA}
\email{Krawzik@shsu.edu}

\author{Colin M. Lawson}
\address{Colin M. Lawson, Department of Mathematics and Statistics \\Stephen F. Austin State University\\
Nacogdoches, Texas, USA}
\email{Colin.Lawson@sfasu.edu}

\author{Christine\ Uhl}
\address{Christine\ Uhl, Department of Mathematics, St.\ Bonaventure University,
St\ Bonaventure, NY 14778}
\email{cuhl@sbu.edu}


\maketitle

\section{Introduction}

The skew group algebra $S(V) \rtimes G$ encodes the action of a group $G$ on the polynomial ring $S(V)$. For certain choices of a vector space $V$ and group $G$, graded Hecke algebras (and their quantum versions) can be considered as graded deformations of $S(V) \rtimes G$ (respectively $S_q(V) \rtimes G$). (See, for example, \cite{GrimleyUhl},
\cite{LevandovskyyShepler},
\cite{NaW16}, 
\cite{SW08}, 
\cite{ShroffWi16}.) Rational Cherednik algebras and symplectic reflection algebras can 
also be constructed as deformations of skew group algebras for $G$ a complex reflection group.
(See, for example, \cite{EG},  
\cite{Gordon08},  
\cite{Norton}.)
Lusztig's graded affine Hecke algebras of \cite{Lusztig89} deform the group action on the vector space and leave the multiplication of the polynomials unchanged. We refer to this type of deformation as Lusztig-type. Drinfeld's deformations of \cite{Drinfeld} instead deform the multiplication of polynomials and leave the group action unchanged.  We refer to this type of deformation as Drinfeld-type. (See, for example, \cite{FGK17} or \cite{SheplerUhl}.)
Later developments considered deforming both the group action and the polynomial multiplication simultaneously. (See 
\cite{SW12} or
\cite{SW14}, 
for example.) Ram and Shepler, in \cite{RamShepler}, showed that in characteristic 0, every deformation of Lusztig-type is isomorphic to a deformation of Drinfeld-type for finite real reflection groups.  Shepler and Witherspoon, in \cite{SW15},  extended this result to the nonmodular setting, where a finite group $G$ has order coprime to the 
characteristic of the field.  It was further shown by Krawzik and Shepler, in \cite{KrawzikShepler}, that in the nonmodular setting, deformations which are of both Lusztig- and Drinfeld-type are isomorphic to a deformation of Drinfeld-type. Shepler and Witherspoon found that in the modular setting, the situation is quite different, and in \cite{SW15} they provide an example of an algebra that is of both Lusztig- and Drinfeld-type but not isomorphic to any algebra of Drinfeld-type.   Shepler and Witherspoon studied the combined Lusztig- and Drinfeld-type deformations of $S \rtimes G$ where $S$ is a Koszul algebra and $G$ is a finite group in \cite{SW18}, using Hochschild cohomology to classify the conditions on the deformations to inform the behavior in the modular and nonmodular settings. 

We extend the study of deformations occurring in the modular setting by specializing to those deformations of $S(V) \rtimes G$ in which $V$ is a vector space over $\mathbb{F}_p$, a finite field of prime order $p\neq 2$, and $G$ is a finite cyclic group of order $p$. In \cite{SW18}, Shepler and Witherspoon provided conditions, in the style of \cite{BG96}, for Poincar\'{e}-Birkhoff-Witt deformations 
when $G$ an arbitrary finite group.  We record these conditions in \cref{SW_conditions}.  These conditions provide an initial refinement of the deformation parameters.  Lawson and Shepler computed Hochschild cohomology of $S(V) \rtimes G$  when $G$ is a cyclic group and char $\FF \neq 2$ in \cite{LS24}. In
\cref{TwistedResolutionSection}, we provide the necessary chain maps to translate the cohomological conditions to be comparable to those of \cite{SW18}.  These conditions alone still leave a lot to be desired.  In \cref{transvection_cohomology}, we further specify these maps for the transvection cyclic group acting on a $2$-dimensional vector space.  
\cref{equations} concludes with a full, albeit difficult to parse, description of the deformations.  In \cref{algebra_equations}, we re-characterize these conditions so that one could choose a single parameter and fully determine a family of deformations.  In \cref{main}, we provide an alternate view of this system, classifying the solution spaces combinatorially and providing a complete formulaic description of all Drinfeld orbifold algebras in our setting.  

\section{Background}
\label{SW_conditions}
Throughout, let $\FF$ be a field of characteristic not $2$ and let $V =\FF^n$ be a finite vector space.  Let $G$ be a finite group.  Unless otherwise noted, $\otimes = \otimes_{\FF}$. 
We will consider actions of $G$ by automorphisms of $V$ which will be denoted, for $g$ in group $G$ and $v$ in vector space $V$, by $\ ^gv$.

\subsection*{Skew group algebras}
For an $\FF$-algebra $S$ with action of a finite group $G$ by automorphisms, we may define the skew group algebra $S \rtimes G=S \otimes \FF G$ as a vector space where for $a, a' \in S$ and $g, g',  \in G$, multiplication is given by $$(a \otimes g)\cdot(a' \otimes g')=a(^g a') \otimes gg'.$$ 
In this paper, we will consider $S=S(V)$ and $G$ a finite cyclic group. The group $G$ will be further specified in \cref{transvection_cohomology}.  When the context is clear, we will suppress the $\otimes$, writing $ag$ rather than $a \otimes g$.

\subsection*{Drinfeld orbifold algebras}
The algebras we explore are described in terms of parameter functions $\kappa: V \times V \rightarrow \FF G \oplus (V \otimes \FF G)$ and $\lambda: \FF G \times V \rightarrow \FF G$ which are bilinear.
We will refer to the projection of $\kappa$ onto the constant part, $\FF G$, as $\kappa^C$ and the projection onto the linear part, $V \otimes \FF G$, as $\kappa^L$.  We define $\mathcal{H}_{\lambda, \kappa}:= (T(V) \rtimes G) /R$, where 
\begin{equation}
\label{DOA}
R = \{gv-{}^gvg -\lambda(g, v), vw-wv-\kappa(v,w) ~|~ g \in G, v, w \in V\}.
\end{equation}
Note that if $\lambda$ and $\kappa$ are trivial for all $g \in G$, $v,w \in V$, we recover the algebra $S(V) \rtimes G$.  We say that $\mathcal{H}_{\lambda, \kappa}$ is a Drinfeld orbifold algebra if it is a Poincar\'e-Birkhoff-Witt (PBW) deformation of $S(V) \rtimes G$. 
The algebra
$\cH_{\ld,\kappa}$ is filtered by degree where each $v \in V$ is assigned degree $1$ and each $g\in G$ is assigned degree $0$.
Then $\cH_{\ld,\kappa}$ 
exhibits the PBW property
when
$$
\gr(\cH_{\ld,\kappa}) \cong S(V)\rtimes G
$$
as graded algebras.
Note that every Drinfeld orbifold algebra 
has an $\FF$-vector space basis given by
$\{v_1^{i_1} v_2^{i_2} \dotsm v_n^{i_n} g:
i_m \in \NN, g\in G\}$
for any $\FF$-vector space basis
$v_1,\ldots, v_n$ of $V$.

\subsection*{PBW Conditions}
Shepler and Witherspoon provided the homological conditions on parameter functions $\lambda$ and $\kappa$ (identified with cochains on a twisted resolution) that yield a Drinfeld orbifold algebra for the more general setting of $G$ a finite group acting linearly on $V$. These conditions, which we refer to as PBW Conditions (1) through (6), we later refine to our specific setting. We give the full statement of their theorem below for completeness. 

\begin{thm}[\cite{SW18} Theorem 6.1]
\label{LiftingConditions}
    Let $G$ be a finite group acting linearly on $V$, a finite dimensional $\FF$-vector space.  Then $\cH_{\ld, \kappa}$ is a PBW deformation of $S(V)\rtimes G$ if and only if 
    \begin{enumerate}
        \item $\displaystyle \ld(gh,v)=\ld(g, \ ^hv)h +g\ld(h,v)$
        
        \item $\displaystyle \kappa^C(\ ^gu,\ ^gv)g-g\kappa^C(u,v)=\ld(\ld(g,v),u)-\ld(\ld(g,u),v)+\sum_{a\in G}\ld(g,\kappa_a^L(u,v))a$
        
        \item $\displaystyle ^g(\kappa^L_{g^{-1}h}(u,v))-\kappa_{hg^{-1}}^L(\ ^gu,\ ^gv)=(\,^hv-\ ^gv)\ld_h(g,u)-(\ ^hu-\ ^gu)\ld_h(g,v),$
        
        \item $\displaystyle 0=2 \sum_{\sigma \in A_3}\kappa_g^C(v_{\sigma(1)}, v_{\sigma(2)})(v_{\sigma(3)}-\ ^gv_{\sigma(3)})+\sum_{a\in G,\sigma\in A_3}\kappa_{ga^{-1}}^L(v_{\sigma(1)}+\ ^av_{\sigma(1)}\kappa_{\alpha}^L(v_{\sigma(2)},v_{\sigma(3)}))-\\
         2\sum_{a\in G,\sigma\in A_3}\kappa_a^L(v_{\sigma(1)},v_{\sigma(2)})\ld_g(a,v_{\sigma(3)}),$ 
        
        \item $\displaystyle 2\sum_{\sigma \in A_3}\ld(\kappa^C(v_{\sigma(1)},v_{\sigma(2)}),v_{\sigma(3)})
        =-\sum_{a\in G,\sigma \in A_3}\kappa^C_{ga^{-1}}
        \left(v_{\sigma(1)}+\ ^av_{\sigma(1)},\, \kappa_a^L(v_{\sigma(2)},v_{\sigma(3)})\right),$
        
        \item $\displaystyle 0=\kappa_g^L(u,v)(w-\ ^gw)+\kappa_g^L(v,w)(u-\ ^gu)+\kappa_g^L(w,u)(v-\ ^gv),$ 
    \end{enumerate}
    
    in $S(V)\rtimes G$, for all $g,h \in G$ and all $u,v,w,v_1,v_2,v_3 \in V$.
\end{thm}

\begin{rem}\label{Rem:dim2}
   A two-dimensional vector space trivially satisfies conditions (4) and (5).
\end{rem}

PBW Conditions (1), (3), and (6) are conditions on the cocyles in Hochschild cohomology, and PBW Conditions (2), (4), and (5) are bracket conditions, which we explain next.

\subsection*{Deformations and Hochschild cohomology}
Let $A$ be an algebra over $\FF.$  Then a deformation of $A$ over $\FF[[t]]$ is an associative algebra with a new multiplication, $*$ , for $a, b \in A$ given by $$a*b = ab + \mu_1 (a \otimes b)t+\mu_2 (a \otimes b)t^2 + ... 
$$
for bilinear maps $\mu_i: A \times A \rightarrow A$ and $i \in \mathbb{N}$.

Hochschild cohomology of an $\FF$-algebra $A$, denoted $\HH^n(A) := \HH^n(A,A)$ is given by $$ \HH^n(A):=\Ext^n_{A^e}(A, A)$$ where $A^e=A \otimes A^{op}$ and $A^{op}$ is the opposite algebra with multiplication given by $a\cdot_{op} b = ba$.  
While the Hochschild cohomology computations benefit from a resolution-independent definition, we are more limited with respect to deformations; each $\mu_i$ arises as a Hochschild 2-cocycle expressed in terms of the bar resolution.  
The  bar resolution of $A$ is a free resolution of $A$ as an $A$-bimodule with 
$$
B(A)_{\DOT}:\quad  ... \xrightarrow{\delta_3} A \otimes A \otimes A \otimes A \xrightarrow{\delta_2} A \otimes A \otimes A \xrightarrow{\delta_1} A \otimes A \xrightarrow{\delta_0} A \rightarrow 0
\,,
$$ 
where the differentials $\delta_m:A \ot A^{\ot m} \ot \to A \ot A^{\ot m-1} \ot A$ are defined by
\begin{equation}
\label{BarDifferentials}
\delta_m(a_0 \otimes a_1 \otimes ... \otimes a_{m+1}) = \sum_{i=0}^m (-1)^{i} a_0 \otimes ... \otimes a_{i}a_{i+1} \otimes ... \otimes a_{m+1}.
\end{equation} 
In the reduced bar resolution, $\overline{B(A)}_{\DOT}$, we replace $A^{\otimes m+2}$ with $A \otimes (\overline{A})^{\otimes m} \otimes A$ where $\overline{A}:= A/(\FF \cdot 1_A)$.
For the reduced bar, the differentials determined by $\delta$ are passed through the quotients, i.e., the image under $\delta_m$ is zero whenever any of the inner tensor components lie in $\FF$.

The associativity of the deformation also requires the following relations on the Gerstenhaber brackets of Hochschild cohomology $$\delta^*_3(\mu_2)=\frac{1}{2}[\mu_1, \mu_1] \textrm{ and } \delta^*_3(\mu_3)=[\mu_1, \mu_2]
\,,$$
where $\delta^*_3$ is the induced differential from $\delta_3$ on the $\Hom(\overline{B(A)}_{\DOT}, A)$ complex.  See \cite{GS86}, \cite{Shakalli12}, and \cite{W19} for more details on the connection between deformations and cohomology.

\subsection*{Graded deformations} 
In general, if $A$ is a graded $\FF$-algebra, a graded deformation of $A$ is one in which the degree of the indeterminate $t$ is set to $1$ and where each multiplication map $\mu_i$ is homogeneous of degree $-i$. 
In this case, we have two gradings on cohomology: one corresponding to the degree of the graded map within a cohomology class and the other corresponding to the degree of the cohomology  class.
The space of maps that are of homological degree $n$ and of graded degree $m$ will be denoted $\HH_{m}^n(A).$

To recall the main result in \cite{LS24}, we first define the map $\displaystyle T: V \rightarrow V$ by $v \mapsto \sum_{h \in G} {}^h v,$ 
and for $h \in G,$ define $$V^h=\textrm{Ker}(1-h)=\{v \in V ~|~ {}^h v = v\} \quad \textrm{ and}\quad 
 V_h=\textrm{Im}(1-h)=\{v-{}^h v ~|~ v \in V\}.$$  
 When $G$ is abelian, we can define a linear character, the determinant of $G$ acting on $V/V^h,$ $$\chi_h: G \rightarrow \FF^{\times} \quad
 \text{by}
 \quad
 \chi_h(g):=\textrm{det}[g]_{V/V^h}
 \quad\text{for all $g \in G$}
 \,.$$

\begin{thm}[\cite{LS24} Theorem 8.1]
\label{MainThm}
Let $G\subseteq \text{GL}(V)$ be a finite cyclic group 
with $V = \FF^n$.
The space of 
infinitesimal graded deformations of $A=S(V)\rtimes G$
is isomorphic as an $\FF$-vector space to
{
$$
\begin{aligned}
\HH_{-1}^2(A)
\cong
(V^G/\im T)^*
\oplus
\big( 
V \ot\, \Wedge^2 V^*
\big)^G
\   \oplus\hspace{-2ex}
\bigoplus_{\substack{h \in G \\ \codim V^h =1 }}
\hspace{-2ex}
\left(\FF \oplus
\big(V/V_h \ot (V^h)^*\big)
\right)^{\chi_h}
\ \oplus\hspace{-2ex}
\bigoplus_{\substack{h\in G \\ \codim V^h =2 }}
\hspace{-2ex}
\left(
V/V_h
\right)^{\chi_h}
.
\end{aligned}
$$
}
\end{thm}

The above theorem identifies cochains that are potential candidates for the first multiplication maps of formal deformations of $S(V)\rtimes G$ when $G$ is a finite cyclic group. However, the cohomology above is computed on a particular resolution, while the lifting conditions of \cref{LiftingConditions} are given on a different resolution.
We take these eligible cochains and define the necessary chain maps to translate the cochains into maps on the resolution of \cite{SW15}.

\section{Resolutions and chain maps for the group}
From here on out, we assume that $G$ is a finite cyclic group acting linearly on a finite dimensional vector space $V$. 
We recall two well-known projective $\FF G$-bimodule resolutions for group algebra $\FF G$: the reduced bar resolution $\overline{B(G)}_{\DOT}$ and a periodic resolution $P_{\DOT}$. 
In addition, we describe an explicit choice of chain maps back and forth between the two resolutions.
These resolutions and chain maps will be used in the next section to create  twisted resolutions for $A=S(V) \rtimes G$, and the chain maps we describe here will be used in creating chain maps between twisted resolutions.

\subsection*{Reduced bar for the group}
The bar resolution is the complex of $F G$-bimodules given by 
\begin{equation}
\label{BarForGroup}
\overline{B(G)}_{\DOT}:
\dots 
\longrightarrow 
\FF G \ot (\overline{\FF G})^{\ot 2} \ot \FF G
\longrightarrow
\FF G \ot \overline{\FF G} \ot \FF G
\longrightarrow \FF G \ot \FF G \longrightarrow \FF G \longrightarrow 0
\,,
\end{equation}
with differentials identical to those given in \cref{BarDifferentials}, and where, as before, $\overline{\FF G}:=\FF G/(\FF \cdot 1_{\FF G}$). 
Note that $G$ need not be cyclic to construct this complex. The canonical $G$-grading on each term $\left(\overline{\FF G}\right)^{(\ot m+2)}$ in the reduced bar is given by the sum of the degrees in each tensor component.

\subsection*{Periodic resolution for the group}
Choose a generator $g$ of the cyclic group $G$. Then the periodic resolution of interest here is the complex of $\FF G$ bimodules given by
\begin{equation}
\label{PeriodicForGroup}
P_{\DOT}:
\dots 
\xlongrightarrow{\gamma} 
\FF G \ot \FF G
\xlongrightarrow{\eta} 
\FF G \ot \FF G 
\xlongrightarrow{\gamma}
\FF G \ot \FF G
\xlongrightarrow{m} 
\FF G 
\longrightarrow 0
\,,
\end{equation}
where $m$ is multiplication, $\gamma=g \ot 1-1 \ot g$, and $\eta = 1\otimes g^{p-1} + g \ot g^{p-2}+ \cdots + g^{p-1} \ot 1$.
Note that this resolution is specific to $G$ a finite cyclic group.
We use the following $G$-grading on $\FF G \ot \FF G$:
For any $h$ in $G$ and $P_i = \FF G \ot \FF G$
for $i \ge 0$, set
\begin{equation}
\label{PeriodicGrading}
(P_i)_{h}=
\begin{cases}
\Span_{\FF}\{a \ot b \ :\ ab = h \}
& \text{ if $i$ is even }
\\
\Span_{\FF}\{a \ot b \ :\ ab = hg^{-1} \}
& \text{ if $i$ is odd}
\,.
\end{cases}
\end{equation}
Such a $G$-grading satisfies the compatibility requirements for twisting resolutions, which we will discuss in the next section.

In the next subsection, we give an explicit choice of chain maps 
between the two resolutions above. 
We identify $G$ with the quotient ring $\FF[x]/(x^{|G|}-1)$ and translate the chain maps described in \cite{GGRSV91} into the context of this paper.

\subsection*{Chain map from the reduced bar to the periodic}
Define the 
map
$
\pi_G: \overline{B(G)}_{\DOT} \longrightarrow P_{\DOT}
$
by
\begin{equation}
\label{PiMapGroup}
(\pi_G)_{n}(1 \ot g^{i_1} \ot \cdots \ot g^{i_{n}} \ot 1)
=
\begin{cases}
\displaystyle\prod_{j=1}^k
\left(1 \ot g^{i_{2j-1}+i_{2j} - |G|}\right)
&
\text{ if }
n=2k
\vspace{2ex}
\\
\left(
\displaystyle\sum_{\ell = 0}^{i_1-1}
g^{\ell} \ot g^{i_1-\ell -1}
\right)
\cdot
\displaystyle\prod_{j=1}^{k}
\left(
1 \ot g^{i_{2j}+i_{2j+1}-|G|}
\right)
&
\text{ if } n=2k+1
\end{cases}
\end{equation}
where $g^{\ell}:=0$ whenever $\ell<0$.
By \cite[Proposition 1.5]{GGRSV91}, this defines a chain map.

\subsection*{Chain map from the periodic to the reduced bar}
In \cite[Proposition 1.5]{GGRSV91}, a chain map in the reverse direction is also given. Define
$
\iota_G: P_{\DOT} \longrightarrow \overline{B(G)}_{\DOT}
$ by
\begin{equation}
\label{IotaMapGroup}
(\iota_G)_n(1\otimes 1) = 
\begin{cases}
\displaystyle\sum_{i_1, \dots, i_k = 1}^{|G|-1}
1 \ot g^{i_k} \ot g
\ot g^{i_{k-1}} \ot g
\ot 
\cdots
\ot g \ot g^{i_1}
\ot g
\ot g^{(k|G| - (\sum_{j=1}^{k} i_j) - k)} 
 &\text{if } n=2k
\vspace{2ex}
\\ 
\displaystyle\sum_{i_1, \dots, i_k =1}^{|G|-1}
1 \ot g \ot g^{i_k} \ot g
\ot g^{i_{k-1}}
\ot g \ot 
\cdots \ot
g^{i_1} \ot g
\ot g^{(k|G| - (\sum_{j=1}^{k} i_j) - k)} & \text{if } n=2k+1.
\end{cases}
\end{equation}

\begin{remark} \em
We note that a few straightforward computations show that the maps defined in \cref{PiMapGroup,IotaMapGroup} are $G$-graded maps of degree zero, and that the composition of chain maps $(\pi_G)_n (\iota_G)_n$ is the identity map on $(\FF G \ot \FF G)_n$ for all $n \ge 0$. 
\end{remark}
Next, we use the above resolutions to create twisted resolutions for $A=S(V)\rtimes G$ and we use the above chain maps to build chain maps between these twisted resolutions.

\section{Twisted product resolutions and chain maps}
\label{TwistedResolutionSection}
In 
\cite{SW14,SW15,SW18, SW19}, 
techniques were developed to create twisted tensor product resolutions for $A=S(V) \rtimes G$. Here, we use two resolutions: one constructed by twisting the reduced bar construction of $\FF G$ (see \cref{BarForGroup}) together with the Koszul resolution $K_{\DOT}$ of $S(V)$, $Y_{\DOT}=\overline{B(G)}_{\DOT} \ot^G K_{\DOT}$, and the other  constructed by twisting the periodic resolution of $\FF G$ (see \cref{PeriodicForGroup}) together with the Koszul resolution of $S(V)$, $X_{\DOT}=P_{\DOT} \ot^G K_{\DOT}$.  
We provide the diagram below to help the reader visualize the resolutions and the maps between them. 
\[
\begin{tikzcd}
X_{\DOT} = P_{\DOT} \ot^G K_{\DOT}:
\ \ \quad 
&[-5ex] 
\cdots 
\to 
&[-6ex] 
\arrow[d,  swap, "\iota",
shift right = .75ex
] 
X_3 
\arrow[r]
&[3ex]
\arrow[d, swap, "\iota",
shift right =2.5ex] 
X_2 
\to 
&[-7ex] 
\cdots 
&[0ex] 
\text{(Hochschild cohomology of $A$ in \cite{LS24})}
\\
Y_{\DOT} = \overline{B(G)}_{\DOT} \ot^G K_{\DOT}: 
&[-5ex]
\cdots 
\to 
&[-6ex] 
\arrow[u, swap, "\pi",
shift right=.75ex] 
\arrow[d, swap,
shift right = .75ex] 
Y_3 
\arrow[r] 
& 
\arrow[u, swap, "\pi",
shift right=-1ex] 
\arrow[d, swap, 
shift right = 2.5ex] 
Y_2  \to &[-7ex] 
\cdots 
&[0ex] 
\text{(lifting conditions in \cref{LiftingConditions})}
\\
Z_{\DOT} = \text{B}(A)_{\DOT}: 
\quad \ \ \ 
&[-5ex]
\cdots 
\to &[-6ex] 
\arrow[u, swap,  
shift right =.75ex] 
Z_3 
\arrow[r] 
& 
\arrow[u, swap,  
shift right=-1ex] 
Z_2  
\to 
&[-7ex] 
\cdots
&[3ex] 
\hspace{0ex}
\text{(for describing deformations)}
\end{tikzcd}
\]

\subsection*{Compatibility conditions for twisting}
We recall here the necessary ingredients needed for twisting resolutions together.
For associative algebras $A$ and $B$ with multiplication maps $m_A: A \otimes A \rightarrow A$ and $m_B: B \otimes B \rightarrow B$, a twisting map 
$\tau: A \otimes B \rightarrow B \otimes A$ is a bijective $\FF$-linear map which commutes with the multiplication maps.  The twisted tensor product algebra $A \otimes_{\tau} B:=A \otimes B$ as a vector space with multiplication determined by $\tau$.  
A skew group algebra is an example of a twisted tensor product algebra where the twisting map $\tau$ is defined by the action of $G$ on $V$. 
For graded algebras, additional constraints require that $\tau$ maintains the grading. Furthermore, with specific compatibility conditions, 
we can combine resolutions of $A$ and $B$ to create a resolution for $A \otimes_{\tau} B$. 
By 
\cite[Proposition 2.20]{SW19}
the Koszul and reduced bar resolution satisfies the compatibility requirements for twisting.  By \cite[Section 4]{LS24}, with the specified grading, the periodic resolution also satisfies the compatibility requirements.

See \cite{BW22}  
for further development of Hochschild cohomology of twisted tensor products and see \cite{KMOOW20} for the Gerstenhaber brackets.

\subsection*{Periodic-twisted-Koszul resolution}
For $A=S(V) \rtimes G$ with $G$ a finite cyclic group acting on $V = \FF^n$, the recipe in \cite{SW14} for twisting the periodic complex in \cref{PeriodicForGroup} with the Koszul complex gives a projective resolution of $A$, which is the total complex
$
X_{\DOT}:=P_{\DOT} \ot^G K_{\DOT}
$.
The resolution, according to \cite[Section 4]{LS24},  is given by
$$
X_n = \bigoplus_{i + j =n} P_i \otimes K_j = (\FF G \otimes \FF G)_i \otimes \left(S(V) \otimes \Wedge^j V \otimes S(V)\right) \cong \left(A \ot \Wedge^j V \ot A\right)
$$ 
where the identification on the right is an $A$-bimodule isomorphism (as in \cite[Section 3]{SW18}) given by
\begin{equation}
\label{PGIsom}
\phi: (a' \ot a) \ot (r \ot y \ot s) \mapsto 
(\,^h r \ot a') \ot (\,^a y) \ot (\,^a s \ot  a)
\end{equation}
for all $(a' \ot a)$ in $(P_i)_h$ (the $h$-th graded component), $r \ot y \ot s$ in $K_j$, and $r,s$ in $S(V)$.
See \cite[Section 4]{LS24} for the differentials on $X_{\DOT}$.

\subsection*{Reduced bar-twisted-Koszul resolution}
The reduced bar-twisted-Koszul resolution, $Y_{\DOT} = \overline{B(G)}_{\DOT} \ot^G K_{\DOT}$ is the total complex 
$$Y_n = \bigoplus_{i+j=n} (\overline{B(G)} \ot^G K)_{i,j} \cong \bigoplus_{i+j=n} \left(A \ot (\overline{\FF G})^{\ot i} \ot \Wedge^j V \ot A\right)
\,.$$
The identification on the right is an $A$-bimodule isomorphism given by 
\begin{equation}
\label{BGIsom}
\psi: (a' \ot x \ot a) \ot (r \ot y \ot s) \mapsto 
(\,^h r \ot a') \ot x \ot (\,^a y) \ot (\,^a s \ot  a)
\end{equation}
for all $(a' \ot x \ot a) \in \left(\overline{B(G)}_i\right)_h$, $r \ot y \ot s \in K_j$, and $r,s \in S(V)$. See \cite[Section 7]{SW15} for details.

\subsection*{Chain maps between twisted resolutions}
Here, we describe a choice of chain maps back and forth between the resolutions $X_{\DOT}$ and $Y_{\DOT}$, allowing the transfer of necessary cohomological information between the resolutions.
Chain map between $Y_{\DOT}$ and $Z_{\DOT}$ are given in \cite[Lemma 4.4]{SW18}, which were used to provide the PBW conditions in \cref{LiftingConditions} on $Y_{\DOT}$.

\begin{lemma}\label{ChainMapLemma}
    There exists chain maps $\iota_{\DOT}: X_{\DOT} \longrightarrow Y_{\DOT}$ and $\pi_{\DOT}: Y_{\DOT} \longrightarrow X_{\DOT}$ of degree zero as graded maps with $\pi_n \iota_n=\id_X$ for all $n\ge 0$. Specifically, let $\iota_n:=\sum_{i+j = n} \iota_{i,j}$ and $\pi_{n}=\sum_{i+j=n} \pi_{i,j}$ be defined
    for all $y$ in $\Wedge^j V$, by
$$
\iota_{i,j}(1 \ot y \ot 1) =
\begin{cases}
    \displaystyle\sum_{\ \ \ i_1, \dots, i_k = 1}^{|G|-1}
1 \ot g^{i_k} \ot g \ot \cdots \ot g^{i_1} \ot g \ot 
\,^{\bar{\beta}} y \ot \bar{\beta}
& 
\text{ if $i=2k$} 
\vspace{1ex}
\\
\displaystyle\sum_{\ \ \ i_1, \dots, i_k = 1}^{|G|-1}
1 \ot g \ot g^{i_k} \ot g \ot \cdots  
 \ot g^{i_1} \ot g \ot 
\,^{\bar{\beta}} y \ot \bar{\beta}
&
\text{if $i=2k+1$}
\,
\end{cases}$$
where $\bar{\beta} := g^{(k|G| - (\sum_{j=1}^k i_j) -k)}$
in $G$,
and for all $y \in \Wedge^j V$ and $g^{i_1},\dots,g^{i_n} \in G$,
$$
\pi_{i,j}(1 \ot g^{i_1} \ot \cdots \ot g^{i_i} \ot y \ot 1) 
= 
\begin{cases}
1 \ot \,^{\beta} y \ot \beta
&
\text{if $i=2k$, where $\beta = \prod_{j=1}^k g^{i_{2j-1}+ i_{2j}-|G|}$}
\vspace{1ex}
\\
\displaystyle\sum_{\ell = 1}^{i_1-1}
\left(
g^{i_1-1-\ell} \ot \,^{g^{\ell} \beta'} y \ot g^{\ell} \beta'
\right)
&
\text{if $i=2k+1$, where $\beta' = \prod_{j=1}^{k} g^{i_{2j}+i_{2j+1}-|G|}$}
\end{cases}
$$
where $g^{\ell}=0$ if $\ell < 0$. 
Then $\iota_n$ and $\pi_n$ are chain maps of degree zero with $\pi_n\iota_n = \text{id}_X$ for all $n \ge 0$.

\end{lemma}
\begin{proof}
Let $\iota_S$ and $\pi_S$ be the identity chain maps back and forth between the Koszul resolution of $S=S(V)$ and itself. 
 Consider the tensor product of maps 
 $$
 \iota_G \ot \iota_S : (P \ot^G K)_{\DOT} \longrightarrow (\overline{B(G)} \ot^G K)_{\DOT}
\quad
\text{and}
\quad
\pi_G \ot \pi_S: (\overline{B(G)} \ot^G K)_{\DOT} \longrightarrow (P \ot^G K)_{\DOT}
\,,$$
where $\iota_G:  P_{\DOT} \longrightarrow \overline{B(G)}_{\DOT}$ and $\pi_G: \overline{B(G)}_{\DOT} \longrightarrow P_{\DOT}$ are as defined in \cref{IotaMapGroup,PiMapGroup}.
Then, using the identifications $\phi$ and $\psi$ in \cref{PGIsom} and \cref{BGIsom} applied to $\iota_G \ot \iota_S$ and $\pi_G \ot \pi_S$, we can define
$$
\iota_{i,j} := \psi \circ  (\iota_G \ot \iota_S) \circ \phi^{-1}: A \ot \Wedge^j V \ot A \longrightarrow A \ot (\overline{\FF G})^{\ot i} \ot \Wedge^j V \ot A
\qquad
\text{for $i, j \ge 0$}
$$ 
and
    $$
    \pi_{i,j}:=\phi \circ(\pi_G \ot \pi_S)\circ\psi^{-1}: A \ot (\overline{\FF G})^{\otimes i} \ot \Wedge^j V \ot A
    \longrightarrow 
    A \ot \Wedge^j V \ot A
    \qquad
    \text{for $i,j \ge 0. $}
    $$
We show that our constructed maps, $\iota_{\DOT}$ and $\pi_{\DOT}$ are chain maps.  
Notice that the maps $d_X$ and $d_Y$ are tensor products of differentials, under the same identifications as above.  That is, $$d_X=\phi(d_P \otimes \id_K + (-1)^i \id_P \otimes d_K)\phi^{-1} \textrm{ and } d_y=\psi(d_B \otimes \id_K +(-1)^i \id_B \otimes d_K) \psi^{-1}.$$  Thus, we can verify that $\iota d_X=d_Y \iota$ precisely because $\iota$ was constructed from chain maps.  Namely
\begin{align*}
    \iota d_X=\psi(\iota_G \otimes \iota_S)\phi^{-1}\phi(d_P \otimes \id_K + (-1)^i \id_P \otimes d_K) \phi^{-1}
    &=
    \psi(\iota_G \otimes \iota_S)(d_P \otimes \id_K + (-1)^i \id_P \otimes d_K) \phi^{-1}
    \\
    &=\psi(\iota_G d_P \otimes \iota_S \id_K + (-1)^i \iota_G \id_P \otimes \iota_S d_K)\phi^{-1}.
\end{align*}
Because $\iota_S$ (the identity chain map) and $\iota_G$ are chain maps, $d_B\iota_G=\iota_G d_P$ and $\iota_S d_K=d_K \iota_S$.  We also know that $\iota_S \id_K=\iota_S=\id_K \iota_S$ and $\iota_G \id_P=\iota_G = \id_G \iota_G$.  Therefore 
\begin{align*}
    \iota d_X
    =
    \psi(d_B\iota_G \otimes \id_K \iota_S  
    &
    + (-1)^i \id_G \iota_G  \otimes d_K \iota_S)\phi^{-1}
    =
    \psi(d_B \otimes \id_K + (-1)^i \id_B \otimes d_K)(\iota_G \otimes \iota_S) \phi^{-1}
    \\
    &
    \qquad
    \qquad
    \qquad
    =\psi(d_B \otimes \id_K + (-1)^i \id_B \otimes d_K)\psi^{-1}\psi(\iota_G \otimes \iota_S) \phi^{-1}
    =d_Y \iota.
\end{align*}
By a comparable computation, again leveraging that $\pi$ is constructed in a similar manner from chain maps, we can see that $\pi$, too, is a chain map.
Lastly, $\pi \iota=\id_X$ since $(\phi \circ(\pi_G \ot \pi_S)\circ\psi^{-1}) \circ (\psi \circ  (\iota_G \ot \iota_S) \circ \phi^{-1}) =\id_X$ as $\pi_G \iota_G=\id_P$ and $\pi_S \iota_S = \id_S$.
\end{proof}

To clarify any ambiguity, we provide the following table of values for chain maps in low degrees.  These will also be useful in the following section.

\begin{table}[H]
\centering
\caption{Chain maps in low degrees   }
\label{tab:lowdimmaps}
$\begin{tabular}{|c|c|c|}
\hline 
   For $n=0$ & &\\
   &$\iota_0$ is the identity map & $\pi_0$ is the identity map\\
\hline
   For  $n=1$  & &  \\
  &$\iota_1(1 \otimes 1)=   1 \otimes g \otimes 1$ & $\pi_1(1 \otimes g^s \otimes 1)=\sum_{\ell=0}^{s-1} g^{s-1-\ell} \otimes g^{\ell}$\\
  &$\iota_1(1 \otimes w_1 \otimes 1)= 1 \otimes w_1 \otimes 1 $& $\pi_1(1 \otimes w_1 \otimes 1)=1 \otimes w_1 \otimes 1$\\
  \hline
  For  $n=2$ &&\\
  &$\iota_2(1 \otimes 1)=\sum_{i=1}^{|G|-1} 1 \otimes g^i \otimes g \otimes 1$&$ \pi_2(1 \otimes g^s \otimes g^r \otimes 1)=1 \otimes g^{s+r-|G|}$\\
  &$\iota_2(1 \otimes w_1 \otimes 1)=1 \otimes g \otimes w_1 \otimes 1$ &$\pi_2(1 \otimes g^s \otimes w_1 \otimes 1)=\sum_{\ell=0}^{s-1} g^{s-1-\ell} \otimes {}^{g^{\ell}}w_1 \otimes g^{\ell}$\\
  &$\iota_2(1 \otimes w_1 \wedge w_2 \otimes 1)=1 \otimes w_1 \wedge w_2 \otimes 1$&
$\pi_2(1 \otimes w_1 \wedge w_2 \otimes 1)=1 \otimes w_1 \wedge w_2 \otimes 1$\\
&&where, in the image of $\pi$, $g^{\ell}=0$ if $\ell<0$.\\
\hline
  
\end{tabular}$
\end{table}

\section{Deformation cohomology for cyclic transvection groups}
\label{transvection_cohomology}
Throughout the remainder of this paper we consider the case where $\FF = \FF_p$ and $G$ is the cyclic group of order $p$ generated by the following nondiagonalizable
reflection 
$$
g = 
\left(
\begin{matrix}
    1 & 1 \\
    0 & 1
\end{matrix}
\right)
\,,
$$
acting on $V=\FF_p^2$ with basis $v_1$ and $v_2$. Notice that $
\,^gv_1 = v_1
\text{ and }
\,^gv_2 = v_1 + v_2
\,.$
The PBW Conditions in \cref{LiftingConditions} are expressed in terms of the Hochschild cohomology of the $Y_{\DOT}$ resolution.  Deformations are generally defined on the bar resolution of the algebra. 
However, as computations on the bar resolution are unwieldy, the more convenient resolution $X_{\DOT}$ defined in \cref{TwistedResolutionSection} is used for specific calculations.
Lawson and Shepler \cite{LS24} provided a description of the Hochschild cohomology governing the graded deformations of $A$. 
We state the result for the transvection group $G=\langle g \rangle$ here.

The example in \cite[Section 9]{LS24} applies \cref{MainThm} to $A=S(V)\rtimes G$ in our setting, giving
\begin{equation}
\begin{aligned}
\label{TransvectionCoh}
\HH^2_{-1}(A) 
&\ \cong\ 
\underbrace{
(V^G)^* \oplus 
(V \ot \Wedge^2 V^*)^G }_{\text{contribution of $1_G$}}
\ \oplus\ 
\bigoplus_{
\substack{h \in G\ \\ h \neq 1_G}}
\underbrace{
\big( \FF_p \oplus
(V/V^G \ot (V^G)^*)
\big)^{G}}_{\text{\ \ 
contribution of reflections}}
\\
&\ =\ 
(\FF_p v_1)^* \oplus 
(\FF_p v_1 \ot \FF_p v_1^*\wedge v_2^*)
\ \oplus\ 
\bigoplus_{
\substack{h \in G\ \\ h \neq 1_G}}
\big( \FF_p \oplus
(V/\FF_p v_1 \ot (\FF_p v_1)^*)
\big)
\ \cong\ \FF_p^{2p}
\, .
\end{aligned}
\end{equation}

So far in this paper we have made an explicit choice of chain maps that facilitate the transfer of cohomological information computed in \cite{LS24} using resolution $X_{\DOT}$ to the resolution $Y_{\DOT}$ (described in \cref{TwistedResolutionSection}).  We now use these to translate the $2$-cocycles given above into 
candidates for $\ld$ and $\kappa$ that satisfy PBW Conditions $(1)$, $(3)$, and $(6)$.

\subsection*{Transferring the cohomology}
Here, we use the chain maps defined in
\cref{TwistedResolutionSection}
to convert $2$-cochains on the periodic-twisted-Koszul resolution of $A=S(V) \rtimes G$ into $2$-cochains on the bar-twisted-Koszul resolution. 
We summarize the conversion of $2$-cocycles in the diagram below. 

$$
\begin{tikzcd}[column sep =small, row sep=large]
X_2=
P_2 \ot^G K_2 =
(A \ot A) 
\oplus 
(A \ot V \ot A) 
\oplus 
(A \ot \Wedge^2(V) \ot A) 
\arrow[d, swap, "\iota_2", shift left = 5ex] 
\arrow[r, "\gamma"] 
& 
A 
&[-2ex] 
\gamma 
\arrow[d, mapsto, "\pi_2^*", shift right= .5ex] 
\\
Y_2=\overline{B(G)}_2 \ot^G K_2 =
(A \ot (\overline{\FF_p G})^{\ot 2} \ot A) 
\oplus 
(A \ot \overline{\FF_p G} \ot V \ot A) 
\oplus 
(A \ot \Wedge^2(V) \ot A) 
\arrow[u, swap, "\pi_2", shift right= 7ex] 
\arrow[r] & A  
&[-2ex] \gamma \circ \pi_2  
\end{tikzcd}
$$

Recall that in \cite{LS24} the Hochschild cohomology governing graded deformations of $A$ was computed using resolution $X_{\DOT}$, and so,
as described in \cite[Section 5]{LS24}, the corresponding space of $2$-cochains on $X_{\DOT}$ of graded degree $-1$ is given by
$
\big(
\Hom_{A^e}(X_2, A)
\big)_{-1}
$, which
decomposes in the following way
into a space of maps on $V$
and a space of maps on $\Wedge^2 V$:
\begin{equation}
\label{CochainDecomp}
\begin{aligned}
\Hom_{\FF_p}(V,\FF_p G) 
\oplus
\Hom_{\FF_p}(\Wedge^2 V, \,  V \ot \FF_p G) 
\, ,
\end{aligned}
\end{equation}
since $X_2 =(A \ot A) \oplus (A \ot V \ot A) \oplus (A \ot \Wedge^2 V \ot A)$
and $\FF_p G$ is the degree $0$ component of $A$.
Here, the only $A^e$-homomorphism from
$A\ot A$ to $A$
of degree $-1$ is the zero map.

\begin{lemma}
\label{TranferringCochains}
Let $\gamma = \lambda' \oplus \alpha \in \big(
\Hom_{A^e}(X_2, A)
\big)_{-1}$ with $\lambda' \in \Hom_{\FF_p}(V, \FF_p G)$ and $\alpha \in \Hom_{\FF_p}(\Wedge^2 V, V \ot \FF_p G)$.
Then $\pi_2^*(\gamma)$ is given by
$$
\pi_2^*(\gamma)(g^i \ot g^j) = 0,
\qquad
\pi_2^*(\gamma)(g^i \ot v) = \sum_{\ell =0}^{i-1} \lambda'(\,^{g^{\ell}}v) g^{i-1},
\qquad 
\text{and}
\qquad
\pi_2^*(\gamma)(u \wedge v) = \alpha(u \wedge v)
$$
for all $g^i, g^j$ in $G$ and all $u,v$ in $V$, where $\pi_2^*(\gamma)$ is viewed as a $2$-cochain on $Y_{\DOT}$.
\end{lemma}
\begin{proof}
By the definition of $\pi_2$ in \cref{tab:lowdimmaps}, $\pi_2^*(\gamma)$ is a $2$-cochain on $Y_{\DOT}=\overline{B(G)}_{\DOT} \ot^G K_{\DOT}$ defined by
$$
\pi_2^*(\gamma)(1 \ot g^i \ot g^j \ot 1) = \gamma(\pi_2(1 \ot g^i \ot g^j \ot 1)) = \gamma\left(1 \ot g^{i+j-|G|}\right) = \gamma(1 \ot 1) g^{i+j-|G|} = 0,
$$
$$
\begin{aligned}
\pi_2^*(\gamma)(1 \ot g^i \ot v \ot 1) 
=
\gamma\left(\sum_{\ell=0}^{i-1} g^{i-1-\ell} \ot \,^{g^{\ell}} v \ot g^{\ell}\right)
&
= 
\sum_{\ell=0}^{i-1}
g^{i-1-\ell} \lambda'(\,^{g^{\ell}}v) g^{\ell}
&
=
\sum_{\ell=0}^{i-1}
\lambda'(\,^{g^{\ell}}v) g^{i-1}
\,,
\end{aligned}
$$
and
$$
\pi_2^*(\gamma)(1 \ot v \wedge w \ot 1) = \gamma(\pi_2(1 \ot v \wedge w \ot 1)) = \alpha(v \wedge w)
\,.
$$
Here, we identify the $A^e$-homomorphism $\pi_2^*(\gamma)$ with an $\FF_p$-homomorphism with the same name.
\end{proof}

\subsection*{Candidates for parameter functions $\lambda$ and $\kappa$}

Now we use the cohomology in \cref{TransvectionCoh} and the chain maps described in \cref{ChainMapLemma} to exhibit candidates for parameter functions $\lambda$ and $\kappa$  
that yield Drinfeld orbifold algebras. 
We will transfer the distinguished cocycles serving as cohomology class representatives of $\HH^2_{-1}(A)$ given in \cref{TransvectionCoh} and then we will transfer the space of related coboundaries described in \cite{LS24}. 
We conclude the section with \cref{LambdaKappaCoboundaryLemma} which  decomposes every cocycle into the sum of cohomology representatives, i.e., a distinguished cocycle,  with a coboundary, which will then be used to further analyze the six PBW Conditions.

\subsubsection*{\bf \em Candidates coming from cohomology representatives}
Consider a generic element $\gamma$ from the vector space description of $\HH^2_{-1}(A)$ in 
\cref{TransvectionCoh}:
$$
\gamma = 
b_0 v_1^* + a_0v_1 \ot v_1^* \wedge v_2^*
+
\sum_{\ell =1}^{p-1}
\left(
a_{\ell} + (-\ell b_{\ell} v_1 + \FF_p v_1) 
\ot v_1^*
\right)
\,,
$$
where $a_i, b_i$ are in $\FF_p$ for $i = 0, \dots, p-1$ and where $v_j^*=\delta_{i,j}(v_i)$ 
for $j=1,2$.
Using the isomorphism $\Phi'$ in the proof of 
\cite[Theorem 8.1]{LS24}, 
we identify $\gamma$ with a distinguished $2$-cocycle $\gamma_{X}$ on $X_{\DOT}$:
$$
\gamma_X: X_2 = (V) \oplus (\Wedge^2 V)
\longrightarrow 
(\FF_p G) \oplus (V \ot \FF_p G) \subset
S(V) \ot \FF_p G
$$
which is a graded map in
$(\Hom_{A^e}(X_2,A))_{-1}$ defined by 
\begin{equation}
\label{eqn0}
\gamma_{X}(v_1)= \sum_{\ell = 0}^{p-1} b_{\ell} \, g^{\ell+1},
\qquad
\gamma_X(v_2)= \sum_{\ell = 1}^{p-1} a_{\ell}\, g^{\ell+1},
\qquad
\text{and}
\quad
\gamma_X(v_1 \wedge v_2)= a_0 v_1 + \sum_{\ell = 0}^{p-1} \ell \, b_{\ell}\, v_2 \, g^{\ell}
\,.
\end{equation}
The next lemma  lifts $\gamma_X$ to a particular $2$-cocycle $\gamma = \gamma_Y$ on $Y_{\DOT}$, the bar-twisted-Koszul resolution,
$$
\gamma: Y_2= (\FF_p G \ot \FF_p G) \oplus (\FF_p G \ot V) \oplus (\Wedge^2 V) \longrightarrow S(V) \ot \FF_p G
\,.
$$
\begin{lemma}
\label{CohRepLemma}
    Let $\gamma$ be a cohomology representative of $\HH^2_{-1}(A)$ in \cref{TransvectionCoh} (i.e. a distinguished cocycle on resolution $X_{\DOT}$ as in \cref{eqn0}).
    Then, $\gamma$ is of the form
    $$
    \gamma(g^i \ot g^j) = 0,
    \quad
    \gamma(g^i \ot v_1) = 
    i\, \left(\sum_{\ell=0}^{p-1} b_{\ell} g^{\ell}\right) \, g^i
 \,,
 \quad
 \gamma(g^i \ot v_2) = 
 \binom{i}{2}
\left(\sum_{\ell = 0}^{p-1}
b_{\ell} g^{\ell}\right) g^i
+
i \left(
\sum_{\ell=1}^{p-1} a_{\ell} g^{\ell}\right) g^i
\,,
    $$
and
$$
\gamma(v_1 \wedge v_2) = 
a_0 v_1 + \sum_{\ell = 0}^{p-1} \ell \, b_{\ell}\, v_2 \, g^{\ell}
\,,
$$
for some constants $a_0, \dots, a_{p-1}$ and $b_0, \dots, b_{p-1}$ in $\FF_p$.
\end{lemma}
\begin{proof}
Apply the chain map $\pi_2$ defined in \cref{ChainMapLemma}
from $X_{\DOT}$ to $Y_{\DOT}$ so that $\gamma$ is the $2$-cochain on $Y_{\DOT}$ given by $\pi_2^*(\gamma_X)$ for $\gamma_X$ as defined in \cref{eqn0} to obtain the  desired result.

\end{proof}

To translate this lemma into the context of $\mathcal{H}_{\lambda, \kappa}$, $\lambda=\gamma|_{\FF_p G \ot V}$ and $\kappa^L = \gamma|_{\Wedge^2 V}$ satisfy PBW conditions (1), (3), and (6).  These candidates need not be restricted to cocycles so we include those conditions coming from coboundaries.

\subsubsection*{\bf \em Candidates coming from coboundaries}
Here we transfer the space of coboundaries from resolution $X_{\DOT}$ to cochains on resolution $Y_{\DOT}$. 
The space of $2$-coboundaries of degree $-1$ on $X_{\DOT}$ is given in \cite{LS24}. 
We state the lemma here.
\begin{lemma}[\cite{LS24} Lemma 5.4]
\label{CoboundaryLemma}
    Let $f: V \rightarrow \FF_p G$ be the $\FF_p$-linear function given by $f(v) = \sum_{j=0}^{p-1} f_j(v)g^j$ for all $v$ in $V$, where $f_j: V \rightarrow \FF_p$ for each $0 \le j \le p-1$.
    Then the coboundary $df$ as a cochain on resolution $Y_{\DOT}$ has the following form: For all $0 \le i \le p-1$,
    $$
    \begin{aligned}
        df(g^i \ot g^i) = 
        df(g^i \ot v_1) = 0,
        \ \ 
        df(g^i \ot v_2) = -i \left(
   \sum_{j=0}^{p-1} f_j(v_1)g^j
   \right)g^i,
   \ \ 
   df(v_1 \wedge v_2) = 
   \sum_{j=0}^{p-1} j \ f_j(v_1)\, v_1 \,g^j
   \,.
    \end{aligned}
    $$
\end{lemma}

\begin{cor}
\label{LambdaKappaCoboundaryLemma}
Let $\lambda, {\color{blue}{\lambda_{\text{coboundary}}}}: \FF_p G \times V \longrightarrow \FF_p G$ and $\kappa^L, {\color{blue}{\kappa^L_{\text{coboundary}}}}:  V \times V \longrightarrow V \ot \FF_p G$ be bilinear functions and fix an $\FF_p$-linear function $f: V \longrightarrow \FF_p G$ with $f(v) = \sum_{j=0}^{p-1} f_j(v)g^j$ for all $v$ in $V$, where $f_j(v) \in \FF_p$ for each $0 \le j \le p-1$. 
Then the functions $\lambda + {\color{blue}{{\lambda_{\text{coboundary}}}}}$ and $\kappa^L+{\color{blue}{\kappa^L_{\text{coboundary}}}}$ satisfy PBW Conditions (1), (3), and (6) if and only if, for all $0 \le i \le p-1$,
  $$
\begin{aligned}
    & \lambda(g^i, v_1)
    +{\color{blue}{\lambda_{\text{coboundary}}(g^i, v_1)}}
    = i \, \sum_{j=0}^{p-1} b_j \, g^{i+j}
    +{\color{blue}{0}}
    =i b g^i 
    +{\color{blue}{0}}
    \qquad \text{for $b=\sum_{j=0}^{p-1} b_j g^j \in \FF_p G$.}
    \,
    \\
    &
    \lambda(g^i, v_2)
    +{\color{blue}{\lambda_{\text{coboundary}}(g^i, v_2)}}
    = \binom{i}{2} \sum_{j=0}^{p-1} b_j\, g^{i+j}
    +
   i \, \sum_{j=1}^{p-1} a_j \, g^{i+j} 
   +{\color{blue}{-i \left(
   \sum_{j=0}^{p-1} f_j(v_1)g^j
   \right)g^i}}
   \,
   \text{ and }
   \\
   &
   \kappa^L(v_1,v_2)
   +{\color{blue}{\kappa^L_{\text{coboundary}}(v_1, v_2)}}
   =
   a_0 v_1 + \sum_{j=0}^{p-1} j \, b_j \, v_2 \, g^{j}
   +{\color{blue}{
   \sum_{j=0}^{p-1} j \ f_j(v_1)\, v_1 \,g^j}}
   \,,
\end{aligned}
$$
for some constants $a_0, \dots, a_{p-1}, b_0, \dots, b_{p-1} $ in $ \FF_p$. 
\end{cor}
%
%

\section{Further restricting candidates for parameter functions}\label{equations}

The goal of this section is to further characterize the parameter functions $\ld$ and $\kappa$ that yield Drinfeld orbifold algebras. We know from \cref{Rem:dim2} that PBW Conditions (4) and (5) are satisfied, leaving  Condition (2) to be examined, for which we begin with the left-hand side and its impact on $\kappa^C$,  That is, we analyze 
$$
\kappa^C(\,^gu, \,^gv)g - g\kappa^C(u,v)
\qquad
\text{for all $g \in G$ and $u,v \in V$}
\,.
$$
\begin{remark}
\label{Lem:Abelian}
    
    If $G$ is abelian and the dimension of $V$ is two, then for all $g \in G$, $\kappa^C(\,^gu, \,^gv)g - g\kappa^C(u,v) = (\det(g)-1)\kappa^C(u,v)g$. 
\end{remark}
Note that when $\det(g) =1$, as in this setting, the left hand side of PBW Condition (2) is zero. This also implies that $\kappa^C \in \FF_pG$ is not restricted by PBW conditions.  Now we have a corollary that allows us to ignore the coboundaries as we examine this condition. 
\begin{cor}
\label{LambdaKappaLemma}
Let $\lambda, {\color{blue}{\lambda_{\text{coboundary}}}}: \FF_p G \times V \longrightarrow \FF_p G$ and $\kappa^L, {\color{blue}{\kappa^L_{\text{coboundary}}}}:  V \times V \longrightarrow V \ot \FF_p G$ be bilinear functions and assume that the functions
$\lambda + {\color{blue}{{\lambda_{\text{coboundary}}}}}$ and $\kappa^L+{\color{blue}{\kappa^L_{\text{coboundary}}}}$ satisfy PBW Conditions (1), (3), and (6).
Then $\lambda + {\color{blue}{{\lambda_{\text{coboundary}}}}}$ and $\kappa^L+{\color{blue}{\kappa^L_{\text{coboundary}}}}$ satisfy PBW Condition (2) if and only if $\lambda$ and  $\kappa^L$ satisfy PBW Condition (2).
\end{cor}

%
\begin{proof}
    In light of the previous remark, we only need to consider the right-hand side. Set $\lambda' := \lambda+\lambda_{\text{cob}}$ and $\kappa' := \kappa^L + \kappa^L_{\text{cob}}$, where these parameter functions are of the form in \cref{LambdaKappaCoboundaryLemma}. 
    We argue that the right-hand side of PBW Condition (2) applied to $\lambda' \oplus \kappa'$ is equal to the right-hand side of PBW Condition (2) applied to $\lambda \oplus \kappa^L$.  
   
Expanding the right-hand side of PBW Condition (2) and using the linearity of $\ld'$ and $\kappa'$ along with \cref{LambdaKappaCoboundaryLemma}, we get

$$
\begin{aligned}
    &
\lambda'(\lambda'(g^i, v_2), v_1) - \lambda'(\lambda'(g^i, v_1),v_2)
+ \sum_{a \in G} \lambda'(g^i, \kappa'_a(v_1,v_2))a
    \\
    &
    =
\lambda(\lambda(g^i, v_2),v_1) 
+ 
\lambda(\lambda_{\text{cob}}(g^i, v_2), v_1) 
- 
\lambda(\lambda(g^i, v_1),v_2) 
-
\lambda_{\text{cob}}(\lambda(g^i, v_1),v_2)
\\
& \ 
+ 
\sum_{a \in G}
\lambda(g^i, \kappa^L_{a}(v_1,v_2))a
+
\sum_{a \in G}
\lambda(g^i, \kappa^L_{\text{cob},a}(v_1,v_2))a
+ 
\sum_{a \in G}
\lambda_{\text{cob}}(g^i, \kappa^L_{a}(v_1,v_2))a
\,.
\end{aligned}
$$
Notice that this is the right-hand side of PBW Condition (2) with four additional terms that sum to zero.
Indeed, by the parameter functions' definitions, according to \cref{LambdaKappaCoboundaryLemma}, 
$$
\begin{aligned}
&
\lambda(\lambda_{\text{cob}}(g^i, v_2), v_1)   
-
\lambda_{\text{cob}}(\lambda(g^i, v_1),v_2)
+
\sum_{a \in G}
\lambda(g^i, \kappa^L_{\text{cob},a}(v_1,v_2))a
+ 
\sum_{a \in G}
\lambda_{\text{cob}}(g^i, \kappa^L_{a}(v_1,v_2))a
\\
&
=
-i \, \sum_{j=0}^{p-1} 
f_j(v_1) \lambda(g^{i+j}, v_1)
-
i\, \sum_{j =0}^{p-1} b_j \lambda_{\text{cob}}(g^{i+j}, v_2)
+ \sum_{j=0}^{p-1} j f_j(v_1) \lambda(g^i, v_1)g^j
+\sum_{j=0}^{p-1} j b_j \lambda_{\text{cob}}(g^i, v_2)g^j
\\
&
= -i \sum_{j=0}^{p-1} f_j(v_1) (i+j) \sum_{\ell =0}^{p-1} b_{\ell} g^{i+j+\ell}
+ 
i\, \sum_{j=0}^{p-1} b_j(i+j) \sum_{\ell =0}^{p-1} f_{\ell}(v_1)g^{i+j+\ell}
+
\sum_{j=0}^{p-1} j f_j(v_1) i \sum_{\ell =0}^{p-1} b_{\ell} g^{i+j+\ell}\\
& 
\ \ \,-
\sum_{j=0}^{p-1} j b_j i \sum_{\ell =0}^{p-1} f_{\ell}(v_1)g^{i+j+\ell}
\,.
\end{aligned}
$$
Distributing across the $(i+j)$ in the first two terms shows that the $i^2$ terms sum to zero, as well as the four remaining $ij$ terms. 
This completes the proof.

\end{proof}

\begin{remark} \em
{\bf (Candidates for parameter functions)} 
\label{LambdaKappaRemark}
\cref{LambdaKappaLemma} allows us to only consider parameter functions
$\lambda: \FF_p G \times V \to \FF_p G$ and $\kappa^L: V \times V \to \FF_p G \ot V$ of the form
$$
\begin{aligned}
    & \lambda(g^i, v_1)
    = i \, \sum_{j=0}^{p-1} b_j \, g^{i+j}
    =i b g^i 
    \qquad \text{for $b=\sum_{j=0}^{p-1} b_j g^j \in \FF_p G$.}
    \,
    \\
    &
    \lambda(g^i, v_2)
    = \binom{i}{2} \sum_{j=0}^{p-1} b_j\, g^{i+j}
    +
   i \, \sum_{j=1}^{p-1} a_j \, g^{i+j} 
   \,
   \text{ and }
   \\
   &
   \kappa^L(v_1,v_2)
   =
   a_0 v_1 + \sum_{j=0}^{p-1} j \, b_j \, v_2 \, g^{j}
   \,,
\end{aligned}
$$
for some fixed constants $a_0, \dots, a_{p-1}, b_0, \dots, b_{p-1} $ in $ \FF_p$ and for all $0 \le i \le p-1$.
We track the constants $a_0, \dots, a_{p-1}$ by considering them as coefficients for $a \in \FF_pG$, i.e., $a=\sum_{j=0}^{p-1} a_j g^j$.
\end{remark}

The lemma below shows that a cochain, in the form above, satisfies PBW Condition (2)  if and only if the $a_j$'s and $b_j$'s are a solution to a particular system of equations.

\begin{lemma}
\label{system}
   The parameter functions $\ld$ and $\kappa$ of the form described in \cref{LambdaKappaRemark} satisfy PBW Condition (2) up to a coboundary if and only if the $a_j,b_j$s satisfy 

    $$
0= 
a_0b_{\ell}
\ +
\sum_{\substack{0 \le j,k < p 
\\ j+k \equiv \ell \  (\text{mod } p)}}
b_k
\,
\left(
- \binom{j+1}{2} b_j
+ j a_j
\right)
\qquad
\text{for each $\ell = 0 ,1, \dots, p-1$}
\,.
$$
\end{lemma}
\begin{proof}
 We fix a cochain $\lambda \oplus \kappa^L$ with constants $a_j$ and $b_j$ in $\FF_p$ and consider PBW Condition (2). 
By expanding the right-hand side and canceling opposite terms, we get:
$$
\begin{aligned}
0 
&=
\lambda(\lambda(g^i, v_2),v_1)
- \lambda(\lambda(g^i,v_1),v_2)
+
\sum_{a \in G}
\lambda(g^i, \kappa_a^L(v_1,v_2)) \, a\\
&=i \binom{i}{2}
\sum_{j,k =0}^{p-1}
b_j b_k \, g^{i+j+k}
+
\binom{i}{2} \sum_{j,k = 0}^{p-1}
j\, b_jb_k \, g^{i+j+k}
+
i \sum_{j=0}^{p-1} \sum_{k=0}^{p-1} j \, a_j b_k\, g^{i+j+k}
\\
&
\hspace{10ex}
- i \sum_{j,k = 0}^{p-1} \binom{i+j}{2} b_jb_k \, g^{i+j+k}
+ 
a_0 \, i \sum_{j=0}^{p-1} b_j\, g^{i+j}
+
\binom{i}{2}
 \sum_{j,k=0}^{p-1} j\, b_jb_k \, g^{i+j+k}.
 \end{aligned}
$$
Re-indexing, factoring, and a change of variable simplifies the above to:
$$
\begin{aligned}
0
&
=
\sum_{\ell=0}^{p-1}
\Bigg(
\sum_{\substack{0 \le j,k < p 
\\ j+k \equiv \ell \  (\text{mod } p)}}
\left(
\left(
i \binom{i}{2} 
+ 2j \binom{i}{2}
- i\binom{i+j}{2}
\right) b_jb_k
+
ija_jb_k
\right)
+ia_0b_{\ell}
\Bigg)
g^{i+\ell}\\
&= 
\sum_{\ell=0}^{p-1}
\Bigg(
\sum_{\substack{0 \le j,k < p 
\\ j+k \equiv \ell \  (\text{mod } p)}}
\left(
- \binom{j+1}{2} b_jb_k
+ j a_jb_k
\right)
+ a_0b_{\ell}
\Bigg) \, i \, g^{i+\ell}
\,.
\end{aligned}
$$
The independence of the group elements (in $\FF_p G$) yields the desired result.

\end{proof}

The resulting system of equations has a cyclic nature to it which we will examine further in the following section as we search for the solutions.

\section{ Relating system of equations to an algebra homomorphism}\label{algebra_equations}
In \cref{transvection_cohomology} we used cohomology to help identify candidates for potential parameter functions $\ld$ and $\kappa$ which, in \cref{equations}, were evaluated with PBW Condition (2) to yield a system of equations that further restrict our choice of parameters.  
Now we develop an algebraic characterization of the solution space for this nonlinear system of $p$ equations. Our approach relies on a key insight:  If we fix $b\in \FF_p G$, then the system becomes linear! 

\subsection*{Cyclically shifted system of equations}
By \cref{system}, we have the system of equations:
$$
0= 
a_0b_{\ell}
\ +
\sum_{\substack{0 \le j,k < p 
\\ j+k \equiv \ell \  (\text{mod } p)}}
b_k
\,
\left(
- \binom{j+1}{2} b_j
+ j a_j
\right)
\qquad
\text{for each $\ell = 0 ,1, \dots, p-1$}
\,.
$$ 
Note that when $j=0$ the summand is zero.  Thus, observe that by setting  
\begin{equation} 
\label{eqn:c def}
c_m = \begin{cases} - \binom{j+1}{2} b_j + ja_j &\text{ for } 1\leq m<p 
\text{ where } m+j \equiv 0 \  (\text{mod } p) ,\\
a_0  &\text{ for } m=0
\end{cases}
\end{equation}
one can see that the coefficients $c_m$ cyclically permute, as seen below:
\begin{align*}
0
&=b_0c_0 + b_1c_1 + \cdots + b_{p-1}c_{p-1} 
&& 
\hspace{-8ex}
(\ell =0)
\\
0&=b_{0}c_{p-1}+b_1c_0 + b_2c_1 + \cdots + b_{p-1}c_{p-2}
&&
\hspace{-8ex}
(\ell=1)
& 
\\
\vdots
&&
\hspace{-8ex}
\\
0&=b_{0}c_{1}+b_1c_2 + b_2c_3 + \cdots + b_{p-1}c_{0}
&& 
\hspace{-8ex}
(\ell = p-1)
\end{align*} 
We provide the following example for $p=3$ throughout the section to motivate computations. 

\vspace{2ex}
\noindent 
{\bf Running Example.}
 Let $p=3$ and 
     fix parameter functions $\lambda: \FF_p G \times V \to \FF_p G$ and $\kappa^L : V \times V \to \FF_p G \ot V$ that  satisfy the six PBW Conditions.  Then, for $a,b \in \FF_pG$, $\ld$ and $\kappa$ are of the form
    $$
\begin{aligned}
    & 
    \lambda(g^i, v_1) 
    =
    i \, \sum_{j=0}^{2} b_j \, g^{i+j}
    \,,
    \quad
    \lambda(g^i, v_2) 
    = 
    \binom{i}{2} \sum_{j=0}^{2} b_j\, g^{i+j}
    +
   i \, \sum_{j=1}^{2} a_j \, g^{i+j} 
   \,,
   \quad
   \kappa^L(v_1,v_2)
   =
   a_0 v_1 + \sum_{j=0}^{2} j \, b_j \, v_2 \, g^{j}
   \,
\end{aligned}
$$ 
by \cref{LambdaKappaRemark}, and \cref{system} implies that 
\begin{equation*}
\label{eqn2}
\begin{aligned}
b_0 a_0  +b_1 (-a_2)  + b_2(a_1-b_1) 
&=
0
\,,
\\
b_0(a_1-b_1) + b_1a_0+b_2(- a_2) 
&=
0
\,,
\text{ and }
\\
b_0 (-a_2)  +b_1 (a_1-b_1) +b_2 a_0 
&=
0
\,.
\end{aligned}
\end{equation*}
Set $c_0 = a_0$, $c_1 = -a_2$, and $c_2 = -b_1 +a_1$,
and substitute these values into the above system to get
\begin{equation*}
\label{eqn3}
\begin{aligned}
b_0c_0 +b_1c_1 + b_2c_2 
&=0
\\
b_0c_2+ b_1c_0+b_2c_1 
&=0\,, \text{ and }
\\
b_0c_1 + b_1c_2+b_2c_0
&=0
\,.
\end{aligned}
\end{equation*}

\subsection*{A new perspective on our system of linear equations}\label{EquationsToKer}
To solve the system of equations that arises from PBW Condition (2), we switch perspectives to an equivalent problem requiring us to find the kernel of a particular ring homomorphism.
Fix an element $b$ in $\FF_p G$ and define the map
\begin{equation}
\label{VarPhiDef}
\varphi_b: \FF_p G \longrightarrow \FF_p G
\quad
\text{by}
\quad
\varphi_b(c) = b \cdot \sigma(c)
\,,
\end{equation}
where $\sigma(c)$ is the automorphism of $\FF_p G$ given by the map $g \mapsto g^{-1}$.
Then $0=\varphi_b(c)$ is {\em exactly} the system of equations in \cref{system}.

\vspace{2ex}

\noindent
{\bf Running Example.}
For $p=3$, fix $b=b_01_G + b_1g+b_2g^2$ in $\FF_p G$.
Then for any 
$c = c_01_G + c_1g+c_2g^2$ in $\FF_p G$, \cref{VarPhiDef} implies that
$$ 
\begin{aligned}
  \varphi_b(c) 
  = b \cdot \sigma(c) 
  &= 
  (b_0+b_1g+b_2g^2) \cdot (c_0+c_1g^2 + c_2g)
  \\  
  &=
  (b_0c_0 + b_1c_1+ b_2c_2)1_G + (b_0c_2+b_1c_0+b_2c_1)g + (b_0c_1+b_1c_2+ b_2c_0)g^2
  \,.
\end{aligned}
$$
In this example, one can see that coefficients of each power of $g$ match the cyclic nature of the equations in \cref{system}.
For $p>3$, these coefficients do indeed align with the equations in \cref{system} when we set $\varphi_b(c)=0$.

\subsection*{Characterization of $\ker(\varphi_b)$}
Now that we know that solving the system of equations is equivalent to finding the $\ker(\varphi_b)$ one can see the importance of characterizing  
$\ker(\varphi_b)$ for any $b \in \FF_p G$. %
In order to give a description of $\ker \varphi_b$, we provide
a list of well-known elementary facts regarding the augmentation map.
\begin{lemma}
\label{AugFactsLemma}
Let $G$ be a finite cyclic group of order $p$, $p > 2$.  Fix $b \in \FF_p G$ and consider the {\em augmentation map} (an $\FF_p G$-homomorphism):
$$
\epsilon: \FF_p G \longrightarrow \FF_p
\qquad
\text{ given by }
\qquad
\sum_{i} c_i \, g^i \mapsto \sum_{i} c_i
\qquad
\text{for all $\sum_i c_i g^i \in \FF_p G$}
\,.
$$
Then the following three properties hold.
\begin{enumerate}
    \item[(a)] $\Ker(\epsilon) = \langle (g-1) \rangle$ as an ideal in $\FF_p G$.
    
    \item[(b)] Every element of $\FF_p G - \Ker(\epsilon)$ is a unit.

    \item[(c)] For $b \in \FF_p G$, there is a unique $k \in \{0, \dots, p\}$ such that $b = (g-1)^k \tilde{b}$ and $\sum_i \tilde{b}_i \neq 0$.
\end{enumerate}
\end{lemma}
\begin{proof}
    Parts $(a)$ and $(b)$ are straightforward exercises. We provide the proof of part $(c)$ here.
Suppose $b \in \FF_p G$ is arbitrary.  Then, because $\{(g-1)^i ~|~ i \in \{0, ..., p-1\}\}$ is a basis for $\FF_p G$, 
$b = \sum_{i=0}^{p-1} z_i (g-1)^i$ for $z_i$ in $\FF_p$.  If $b \not\equiv 0$, then there is some smallest $k \in \{0,\dots, p-1\}$ such that $z_k \neq 0$. So we can then factor out $(g-1)^k$ and we have
$b = (g-1)^k \, \tilde{b}$ and $\tilde{b}$ has a nonzero constant term $z_k$. 
\cref{AugFactsLemma}(a) implies that $\tilde{b}$ is not in the kernel of $\epsilon$. Thus by \cref{AugFactsLemma}(b) the element $\tilde{b}$ is invertible.
Note that if $b\equiv 0$ and then $b = (g-1)^p \tilde{b}$ holds for any $\tilde{b}$ thus the equation holds for a $\tilde{b}$ which is invertible in $\FF_p G$.
\end{proof}

Now we are ready to state the result characterizing the kernel of the map $\varphi_b$.
\begin{prop}
\label{prop:KerPhi}
    Fix $b \in \FF_p G$. 
    Then
        $$
       \ker(\varphi_b)= 
        \Span_{\FF_p}\left\{
        (g-1)^{p-j}
        \ |\ 0 \le j \le k
        \right\}
        \,,
        $$
        for the unique $k \in \{0,\dots, p\}$ such that $b = (g-1)^k\, \tilde{b}$ with $\sum_i \tilde{b}_i \neq 0$.
\end{prop}
\begin{proof}
    Fix $b \in \FF_p G$. 
    We consider first the case when $b \equiv 0$ (i.e. $k=p$).  Then $\varphi_b \equiv 0$ and so $\ker(\varphi_b) = \FF_p G = \Span_{\FF_p} \{(g-1)^0, (g-1), \cdots ,(g-1)^p\}$\footnote{Note here we've added $0=(g-1)^{p}$ to the basis of $\FF_pG$.}.
    Now consider $b \nequiv 0$ (i.e. $k < p$).
    Then $c \in \ker(\varphi_b)$ if and only if 
    $
    0 = \varphi_b(c) = b\cdot \sigma(c) 
    $. 
    By \cref{AugFactsLemma}(c), there is a unique $k \in \{0, \dots, p-1\}$ such that $b = (g-1)^k \tilde{b}$ and $\tilde{b}$ is invertible,
 so 
    $$
    \begin{aligned}
        0 = b \cdot \sigma(c)
        = (g-1)^k \tilde{b} \cdot \sigma(c)
        \,.
    \end{aligned}
    $$
    Thus, 
    $(g-1)^k \sigma(c)=0$ as $\tilde{b}$ is invertible. 
Since $\sigma$ is an involution, $\sigma^2$ is the identity map and thus
$$
0 = 
\sigma(0) 
= 
\sigma((g-1)^k\, \sigma(c))
=
\sigma((g-1)^k) \cdot \sigma^2(c)
=
\sigma((g-1)^k) \cdot c
\,.
$$
We argue that $\sigma((g-1)^k) \cdot c = (g-1)^k\, \overline{b} \cdot c$ for some $\overline{b} \in \FF_p G$ with $\overline{b}$ invertible.
Observe that 
$$
\sigma(g-1)
=g^{-1}-1 
= 
g^{p-1}-1
=
(g-1)(g^{p-2}+g^{p-3}+\cdots+g+1)
\,,
$$
where $\epsilon(g^{p-2}+g^{p-3} + \cdots + g+1) = p-1 \neq 0$. 
Thus, since $\FF_p G$ is commutative, 
$$
\sigma((g-1)^k)
= (\sigma(g-1))^k
=
((g-1)(g^{p-2}+g^{p-3}+\cdots+g+1))^k 
= 
(g-1)^k \, \overline{b}
\,,
$$
where 
$\overline{b}=(g^{p-2}+g^{p-3}+\cdots+g+1)^k$, which does not lie in $\ker(\epsilon)$ as $\FF_p$ has no zero-divisors. Therefore, $\overline{b}$ is invertible by \cref{AugFactsLemma}(b), and so 
\begin{equation}
\label{Eqn1}
(g-1)^k \cdot c = 0
\,.
\end{equation}
%
Now, as $\FF_p G =\Span_{\FF_p}\{(g-1)^{\ell}\ |\ 0 \le \ell \le p-1 \}$, write $c=\sum_{\ell=0}^{p-1} c'_{\ell}(g-1)^{\ell}$ for some constants $c'_0, \dots, c'_{p-1}$ in $\FF_p$ and observe that since $(g-1)^j = 0$ for $j \ge p$, \cref{Eqn1} implies that
$$
\begin{aligned}
0=(g-1)^k \cdot c 
&
= (g-1)^k\cdot \big(c'_0 + c'_1(g-1)+ c'_2(g-1)^2+ \cdots+ c'_{p-1}(g-1)^{p-1}\big)
\\
&
=
c'_0 (g-1)^k + c'_1(g-1)^{k+1}+ c'_2 (g-1)^{k+2} + \cdots + c'_{p-k-1} (g-1)^{p-1}
\,.
\end{aligned}
$$
So, since $(g-1)^k, \dots, (g-1)^{p-1}$ are (nonzero and) linearly independent, we must have $c'_0 = \cdots = c'_{p-1}=0$ and thus $c$ lies in $\Span_{\FF_p}\left\{
        (g-1)^{p-j}
        \ |\ 0 \le j \le k \right\}$ as claimed.

\end{proof}

Fixing a value of $b$ allows one to quickly generate solutions with this characterization that define Drinfeld orbifold algebras.  We demonstrate this by revisiting our previous example.

\vspace{2ex}

\noindent
{\bf Running Example.}
For $p=3$, fix $b \in \FF_p G$.  Say $b= 1-g$.
Applying \cref{prop:KerPhi},
tells us $$\ker(\varphi_{1-g}) = \textrm{span}_{\FF_p}\{(g-1)^2\} = \textrm{span}_{\FF_p}\{1+g+g^2\}=\{0,\, 1+g+g^2,\, -1-g-g^2\}.$$
For any given $c\in \ker(\varphi_{1-g})$ with $c$ defined by $c_0 = a_0$, $c_1 = -a_2$, and $c_2 = -b_1 +a_1$, as seen earlier in our running example, we get a solution to the system. To provide a specific example, we consider $c=-1-g-g^2$. In this case $a_0=-1$, $a_2 = 1$ and $a_1+1=-1$, i.e., $a_1 = 1$. Thus $a=-1+g+g^2$ is one of the three solutions corresponding to our fixed $b$. Each solution yields a different $\lambda(g^i,v_2)$ and $\kappa^L(v_1,v_2)$ parameter. For our chosen $c$ and $0\leq i\leq 2$, our parameters are:
$$
\lambda(g^i,v_1) = i(1-g)g^i, \quad \lambda(g^i,v_2) = {i \choose 2}(1-g)g^i + i( g^{i+1} +g^{i+2}),\quad  \text{ and }\quad 
\kappa^L(v_1,v_2) =-v_1 -gv_2,
$$
\noindent and, with choice of $\kappa^C(v_1,v_2) \in \FF_p G$ and $f$ from \cref{LambdaKappaCoboundaryLemma}, one can use \cref{DOA} to build $\cH_{\ld,\kappa}$.


Ignoring $\kappa^C$ and coboundaries, we computed the 81 solutions for $a$ and $b$ when $p=3$.  We include the this solution space in \cref{tab:binomials}.  Including $\kappa^C$ and coboundaries, there are $59,049=3^{10}$ solutions.

\begin{center} 

\begin{table}[H]
\centering
\caption{The solution space for $p=3$ }
\label{tab:binomials}

\medskip
\begin{tabular}{|c|c|}
\hline
    $\textbf{\textit{b}} $&$ \textbf{\textit{a}} $ \\
 \hline 
 $0$& $\FF_3G$\\ \hline
 $b_0+b_1g+b_2g^2 $ with  $\sum_i b_i \neq 0$ & $b_1g$ \\
\hline
 $-(g-1)$,~ $g(g-1)$ & $-g$ , $1-g^2$ , $-1+g+g^2 $\\ 
    $g-1$ ,~ $-g(g-1)$ & $g$ , $-1+g^2,$ $1-g-g^2$ \\
    $(g+1)(g-1) $& $0$ , $1+g-g^2$ , $-1-g+g^2$ \\
    $-(g+1)(g-1)$ & 0 , $-1-g+g^2$ , $1+g-g^2$ \\
    \hline
    $(g-1)^2$ & $1, g, -g^2, -1-g, -g+g^2, -1+g^2, -1+g-g^2,$ $1-g-g^2$, $1+g+g^2$ \\
    $-(g-1)^2$ & $-1, -g, g^2, 1+g, g-g^2, 1-g^2, 1-g+g^2, -1+g+g^2, -1-g-g^2$ \\
\hline
\end{tabular}
\end{table}
\end{center}
Notice that in \cref{tab:binomials}, the first line is the single case where \cref{prop:KerPhi} has $k=p$, when $b=0$, and thus all $a \in \FF_3 G$ are solutions.  The next line corresponds to the $18=3^3-3^2$ values of $b$ in which $\sum_i b_i \neq 0$ and thus $k=0$ in \cref{prop:KerPhi}, requiring $c=0$ and $a_1=b_1$.  The four subsequent lines are the $6 = 3^1(3-1)$ values of $b$ in which $k=1$, in which each fixed $b$ yields $3$ solutions.  The last two lines are the $b$ values in which $k=2$, in which each fixed $b$ yields $3^2$ solutions.

\section{Classifying the Drinfeld orbifold algebras}\label{main} 

This section begins with combinatorial results about the dimension of the solution space and is followed up with our main theorem, which classifies the Drinfeld orbifold algebras, $\cH_{\lambda, \kappa}$ for a fixed $b \in \FF_p G$.
The approach is to view the fixed $b$ in terms of its $k$-class, as in \cref{prop:KerPhi}, which turns out to be useful not just for counting the number of solutions but also for counting the $b$'s themselves.

\subsection*{Counting the solution space}
We make use of the class structure implied by \cref{prop:KerPhi} to determine the size of the solution space. 
\begin{prop}\label{count}
    The solution space for $a$ and $b$ as in \cref{system} has cardinality $p^{p+1}$.
\end{prop}

\begin{proof}
For an arbitrary $0 \neq b \in \FF_p G$, \cref{AugFactsLemma}{(c)} states that $b=(g-1)^k\tilde{b}$ for a unique $k \in \{0,...p-1\}$.  
In terms of the $\{g^i ~|~ 0 \leq i \leq p-1\}$-basis, $b$ is a polynomial of degree at most $p-1$. 
Because $b=(g-1)^k\tilde{b}$, $\tilde{b}$ must be a polynomial of degree at most $p-k-1$. 
In order for $\sum_i \tilde{b}_i \neq 0$, of the $p-k$ coefficients, there are there are $p$ choices for $p-k-1$ coefficients and $p-1$ choices for the remaining coefficient (to exclude the one scenario in which the coefficients would sum to 0). Thus the size of the class of $b$'s with $b=(g-1)^k\tilde{b}$ for a fixed $k$ is $p^{p-k-1}\cdot(p-1)$.

\begin{table}[H]
\begin{center}
\caption {Size of $a$ class for each $b$ class}
$\begin{tabular}{|c|c|c|c|c|}\hline
   $\mathbf{\textit{b}}$ \bf{class}  &  $|\mathbf{\textit{b}}$ \bf{class}$|$ & $\mathbf{\textit{a}}$ \bf{class for} $b$ \bf{fixed} & $|\mathbf{\textit{a}}$ \bf{class for }$b$ \bf{fixed} $|$& $|\mathbf{\textit{a}}$  \bf{class for} $b$ \bf{class} $|$\\
 \hline 
 $0 $& $1$ & $\FF_p G$ & $p^p $& $p^p$ \\
 $\sum_i b_i \neq 0$ & $p^{p-1} \cdot (p-1)$ &  & $1$ & $p^p-p^{p-1}$\\
 $(g-1)\tilde{b}$ & $p^{p-2}\cdot (p-1)$ &  & $p$ &$ p^p-p^{p-1}$ \\
 $(g-1)^2\tilde{b}$ & $p^{p-3}\cdot (p-1)$ &  & $p^2$ & $p^p-p^{p-1}$\\
 $(g-1)^3\tilde{b}$ & $p^{p-4}\cdot (p-1)$ &  & $p^3$ &$ p^p-p^{p-1}$ \\
 $(g-1)^4\tilde{b}$ & $p^{p-5}\cdot (p-1) $&  & $p^4 $& $p^p-p^{p-1}$\\
 $\vdots$ & $\vdots$ &$ \vdots$ & $\vdots$ & $\vdots $\\
 $(g-1)^{p-1}\tilde{b}$ & $p-1$ & &$ p^{p-1}$ & $p^p-p^{p-1} $\\
\hline
\end{tabular}$
\end{center}
\end{table}
\noindent For a fixed $b$, the $|a \textrm{ class}|$ is given by $|\ker(\varphi_{b})|=p^k.$ Summing the final column, we find the cardinality of the solution space for the $a$'s and $b$'s as in \cref{system} is $p^{p+1}$.
\end{proof}

\subsection*{Characterizing the solution space}
The careful reader likely noticed the large gap of information in the table for the general case of $\FF_p$.  We use the remainder of this section to carefully construct the missing $a$ solution classification.  As in \cref{algebra_equations}, after fixing  $b$, we find $a$ in two steps:  First, \cref{prop:KerPhi} gives us a description of $c$  which we then unpack using \cref{eqn:c def} to solve for $a$ in terms of $b$.

If $\sum_i b_i \neq 0$, then, in \cref{prop:KerPhi} $k=0$, thus $c=0.$  Therefore, by definition of $c_m$ (see \cref{eqn:c def})
$${j+1 \choose 2}b_{j}=ja_{j}$$ and thus, by Fermat's Little Theorem, $$j^{p-2}{j+1 \choose 2}b_{j}=a_{j}.$$  Therefore, for a fixed $b$ with $\sum_i b_i \neq 0$, we have a unique solution $$a=\sum_{j=1}^{p-1} j^{p-2}{{j+1} \choose 2} b_j g^j=\sum_{j=1}^{p-2} j^{p-2}{{j+1} \choose 2} b_j g^j.$$

We demonstrate one more small case before generalizing to arbitrary $k$.  If $b=(g-1)\tilde{b},$ (that is, $k=1$,) then by \cref{prop:KerPhi}, $c \in \Span_{\FF_p}\{(g-1)^{p-1}\}$.  Therefore $$c=d\left(1-{p-1 \choose 1}g+{p-1 \choose 2}g^2 - {p-1 \choose 3}g^3 + ... + g^{p-1}\right)=d\sum_{j=0}^{p-1} {p-1 \choose j}(-1)^{j}g^{j}$$
for some $d \in \FF_p$.  We also know that, 
for $p-j \in \{1, ..., p-1\}$, by definition of $c_{p-j}$ in \cref{system} $$c_{p-j}=-{j+1 \choose 2}b_{j}+ja_{j}.$$  Comparing $g^{p-j}$ coefficients, we get that $$(-1)^{p-j}{p-1 \choose p-j}d=-{j+1 \choose 2}b_{j}+ja_{j}.$$  Therefore $$a_0=d \textrm{ and } a_j=j^{p-2}\left((-1)^{p-j}{p-1 \choose p-j}d+{j+1 \choose 2}b_j \right) \textrm{ for } j \in \{1,...,p-1\}, d \in \FF_p.$$

To handle the general case, we need to make a change of basis from the $g-1$ basis of \cref{prop:KerPhi} to the $g$ basis of \cref{system} (or, equivalently \cref{eqn:c def}).  By the Binomial Theorem, for $i \in \{1, ..,. p-1\},$  
$$(g-1)^{i}=\sum_{j=0}^i {i \choose j}g^j(-1)^{i-j}=(-1)^{i}\sum_{j=0}^{p-1} {i \choose j}(-1)^{j}g^j$$
where, in the last equality, we leverage that ${i \choose j}=0$ if $i<j$. 
If $b=(g-1)^k\tilde{b}$ such that $\sum_i \tilde{b}_i \neq 0$, then, for $d_1, d_2, d_3, ..., d_k \in \FF_p$, $$a_0=d_1 - d_2 + d_3 -d_4+\cdots +(-1)^{k-1}d_k$$ and, for $j \in \{1,..., p-1\}$,  $$a_j=j^{p-2}\left((-1)^{p-j} \left({p-1 \choose p-j}d_1 - {p-2 \choose p-j}d_2 +  \cdots + (-1)^{k+1}{p-k \choose p-j}d_k\right)+{j+1 \choose 2}b_j\right).$$

\subsection*{Classifying Drinfeld orbifold algebras}
Fix $b \in \FF_p G$ such that $b=(g-1)^k\tilde{b}$ with $\sum_i \tilde{b}_i \neq 0$, $k$ values $d_1, d_2, d_3, ..., d_k \in \FF_p$, and $\kappa^C \in \FF_p G$ arbitrary.  Using this description of $a$, we can leverage \cref{LambdaKappaRemark} to fully describe the corresponding Drinfeld orbifold algebra.  
Before we state the main result, we simplify notation.  Let $\mathbf{d}=(d_1, .., d_k)$ as fixed above. Define the function $$\mu(\mathbf{d}, j)=(-1)^{p-j}\left( {p-1 \choose p-j} d_1 - {p-2 \choose p-j} d_2 + ... + (-1)^{k+1}{p-k \choose p-j}d_k\right)$$ for $j \in \{0,1,...,p-1\}$.


\begin{thm}\label{Thm:GLUK}
Let $p$ be an odd prime,
$G$ the cyclic group of order $p$ generated by $g$, a nondiagonalizable
reflection, 
$$
g = 
\left(
\begin{matrix}
    1 & 1 \\
    0 & 1
\end{matrix}
\right)
\,,
$$
acting linearly on $V=\FF_p^2$ with basis $v_1$ and $v_2$.  There are $p^{2p+1}$ Drinfeld orbifold algebras (up to a coboundary) given by the following procedure.  

Fix $b \in \FF_p G$ such that $b=(g-1)^k\tilde{b}$ with $\sum_i \tilde{b}_i \neq 0$, $d_1, d_2, d_3, ..., d_k \in \FF_p$, and $\kappa^C \in \FF_p G$. Then $\mathcal{H}_{\lambda, \kappa}$ with 

$$\lambda(g^i,v_1)=ibg^i,$$
$$\lambda(g^i, v_2)=\sum_{j=0}^{p-1} \left({i \choose 2} + i j^{p-2}{j+1 \choose 2} \right)b_j g^{i+j} + i\sum_{j=0}^{p-1} j^{p-2} \mu(\mathbf{d}, j) g^{i+j},$$
and $$\kappa(v_1, v_2)=(d_1 - d_2 + ... + (-1)^{k+1} d_k)v_1 + \sum_{j=0}^{p-1}j b_j v_2 g^j + \kappa^C.$$
is a Drinfeld orbifold algebra. 
\begin{proof}
From \cref{count}, we have $p^{p+1}$ solutions for the $a$ and $b$ values.  The value of $\kappa^C$ is allowed to be arbitrary, yielding an additional $p^p$ solutions for each of the previous solutions.
\end{proof}


\end{thm}

\begin{cor} 
    Let $p$ be an odd prime and  $G$ as before.  There are $p^{3p+1}$ Drinfeld orbifold algebras, $\mathcal{H}_{\lambda,\kappa}$, and each is given by the following procedure.  

    Fix $b \in \FF_p G$ such that $b=(g-1)^k\tilde{b}$ with $\sum_i \tilde{b}_i \neq 0$, $d_1, d_2, d_3, ..., d_k \in \FF_p$, an $\FF_p$-linear function $f: V \to \FF_p G$, and $\kappa^C \in \FF_p G$. Set 

$$\lambda(g^i,v_1)=ibg^i,$$
$$\lambda(g^i, v_2)=\sum_{j=0}^{p-1} \left({i \choose 2} + i j^{p-2}{j+1 \choose 2} \right)b_j g^{i+j} + i\sum_{j=0}^{p-1} j^{p-2} \mu(\mathbf{d}, j) g^{i+j}+{\color{blue}{-i \left(
   \sum_{j=0}^{p-1} f_j(v_1)g^j
   \right)g^i}}$$
and $$\kappa(v_1, v_2)=(d_1 - d_2 + ... + (-1)^{k+1} d_k)v_1 + \sum_{j=0}^{p-1}j b_j v_2 g^j + \kappa^C
   +{\color{blue}{
   \sum_{j=0}^{p-1} j \ f_j(v_1)\, v_1 \,g^j}}.$$

\end{cor}

\begin{proof}
By \cref{LambdaKappaLemma}, if we add on the coboundary terms as specified in \cref{LambdaKappaCoboundaryLemma}, $\lambda$ and $\kappa$ will still meet all of the conditions of \cref{LiftingConditions}.  We have an additional $p$ choices for each of the $f_j(v_1)$ which yields an additional $p^p$ Drinfeld orbifold algebras for each of the $p^{2p+1}$ representatives from the previous result.
\end{proof}

One big (but notoriously difficult) question is how many Drinfeld orbifold algebras there are up to isomorphism?
As noted in the introduction, interesting things happen in characteristic $p$.  Perhaps some of the Drinfeld orbifold algebras in this paper will provide additional examples of deformations of Drinfeld- and Lusztig-type that are not isomorphic to one of Drinfeld type.  Currently there is only one known example \cite{SW15}. 



\end{document}